\documentclass[11pt,reqno]{amsart}
\usepackage{amsthm,amssymb,amsmath}

\usepackage[T1]{fontenc}
\usepackage{lmodern}

\newtheorem{theorem}{Theorem}
\newtheorem*{theorem*}{Theorem}

\theoremstyle{definition}

\newtheorem*{definition*}{\sc Definition}

\newtheorem{remark}{\sc Remark}

\newtheorem*{example*}{\sc Example}

\newcommand{\loc}{{\rm loc}}

\begin{document}

\title{Laplacian with singular drift in a critical borderline case}

\author{D.\,Kinzebulatov}

\makeatletter
\@namedef{subjclassname@2020}{%
  \textup{2020} Mathematics Subject Classification}
\makeatother

\begin{abstract}
We establish well-posedness and regularity for parabolic diffusion equation 
in the case when the singularities of a general drift reach the critical magnitude. The latter dictates the need to work in an Orlicz space situated between all $L^p$ and $L^\infty$.
\end{abstract}

\email{damir.kinzebulatov@mat.ulaval.ca}

\address{Universit\'{e} Laval, D\'{e}partement de math\'{e}matiques et de statistique, Qu\'{e}bec, QC, Canada}



\thanks{The research of the author is supported by the NSERC (grant RGPIN-2017-05567).}

\maketitle

\section{Introduction and result}

\subsection{Introduction}
The paper investigates  well-posedness and other properties of the parabolic diffusion equation 
in the case when the magnitude of critical-order singularities of a general drift attains the critical value.
 This  equation arises, in particular, in the study of particle systems
when the strength of singular attracting interactions of Hardy type reaches (and surpasses) the critical threshold. The latter brings qualitative change in the behaviour of the particles, i.e.\,they start to agglomerate and one can prove only local  in time existence of solution for the corresponding stochastic differential equation (SDE). To be more precise, consider SDE
\begin{equation}
\label{hardy_sde}
X_t-x=-\sqrt{\delta}\frac{d-2}{2}\int_0^t|X_s|^{-2}X_sds + \sqrt{2}B_t,
\end{equation}
where $B_t$ is the standard Brownian motion in $\mathbb R^d$, $d \geq 3$,
with initial condition $X_0=x$ and the drift that pushes the trajectories towards the origin.  (This SDE is obtained from a two-particle system with attracting interaction, $X_t$ being the difference between their trajectories, see Remark \ref{mv_rem}.) If $X_0=0$ and $\delta>4(\frac{d}{d-2})^2$, then this SDE does not have a weak solution. If $\delta>4$, then for every $X_0 \neq 0$ $X_t$ arrives at the origin  in finite time with positive probability.
Informally, the attraction to the origin by the drift is too strong. The proof can be found in \cite{BFGM}.

On the other hand, in \cite[Theorem 1.3]{KiS}, Sem\"{e}nov and the author proved the following  result. Consider SDE
\begin{equation}
\label{sde}
X_t-x=-\int_0^tb(X_s)ds + \sqrt{2}B_t, \quad x \in \mathbb R^d, 
\end{equation}
where $B_t$ is the standard Brownian motion in $\mathbb R^d$, $d \geq 3$, and drift $b:\mathbb R^d \rightarrow \mathbb R^d$ is \textit{form-bounded}, i.e.\,$|b| \in L^2_{\loc}$ and the following quadratic form inequality holds:
\begin{equation}
\label{fbb}
\|b\varphi\|_2^2 \leq \delta\|\nabla \varphi\|_2^2+c_\delta\|\varphi\|_2^2 \quad \forall \varphi \in C_c^\infty
\end{equation}
for some constants $\delta>0$ and $c_\delta<\infty$ (here and in what follows, $\|\cdot\|_p:=\|\cdot\|_{L^p}$). We abbreviate \eqref{fbb} as $b \in \mathbf{F}_\delta$. A broad sufficient condition for $b \in \mathbf{F}_\delta$ is a scaling-invariant Morrey class
\begin{equation*}
\|b\|_{M_{2+\varepsilon}}:=\sup_{r>0, x \in \mathbb R^d} r\biggl(\frac{1}{|B_r(x)|}\int_{B_r(x)}|b|^{2+\varepsilon}dx \biggr)^{\frac{1}{2+\varepsilon}}<\infty,
\end{equation*}
for $\varepsilon>0$ fixed arbitrarily small. Here $B_r(x)$ denotes the ball of radius $r$ centered at $x$. Then the form-bound $\delta=c_d\|b\|_{M_{2+\varepsilon}}$ \cite{F}, see also \cite{CF}.
In particular, vector fields $b$ with entries in $L^d$ or in the weak $L^d$ classes are form-bounded\footnote{The former inclusion is easily seen directly: if $|b| \in L^d$, then, for every $\varepsilon>0$, we can represent $b=b_1+b_2$, where $\|b_1\|_d<\varepsilon$ and $\|b_2\|_\infty<\infty$. So, we obtain, using the Sobolev inequality,
\begin{align*}
\|b\varphi\|_2^2 
 \leq 2\|b_1\|_d^2 \|\varphi\|_{\frac{2d}{d-2}}^2 + 2\|b_2\|_\infty^2 \|\varphi\|_2^2 \leq C_S 2 \|b_1\|_d^2 \|\nabla \varphi\|_2^2 + 2\|b_2\|_\infty^2 \|\varphi\|_2^2,
\end{align*}
so $b \in \mathbf{F}_{\delta}$ with $\delta=C_S 2\varepsilon$. Thus, $\delta$ can be chosen arbitrarily small. In this sense, class $|b| \in L^d$ is sub-critical.}. The form-boundedness condition in fact appears already in the Lax-Milgram theorem since it provides coercivity for the corresponding Dirichlet form, see discussion in the end of Section \ref{res_sect}. A model example of a form-bounded drift $b$ with $|b| \not \in L^d$ is 
\begin{equation}
\label{hardy}
b(x)=\pm \sqrt{\delta}\frac{d-2}{2}|x|^{-2}x.
\end{equation}
The fact that this $b \in \mathbf{F}_\delta$ with $c_\delta=0$ is the well known Hardy inequality. So, SDE \eqref{hardy_sde} is a special case of \eqref{sde}. We emphasize that the form-boundedness of $b$ does not require the existence of even locally summable divergence ${\rm div\,}b$ or any kind of symmetry of $b$.
Returning to the result in \cite[Theorem 1.3]{KiS}, there the authors proved that if $$\delta<4,$$ then SDE \eqref{sde} has a martingale solution for every initial point $x \in \mathbb R^d$. The proof in \cite{KiS} used some elements of De Giorgi's method ran in $L^p$, $p>\frac{2}{2-\sqrt{\delta}}$.

The present papers deals with the borderline case $$\delta=4$$ at the level of the corresponding to \eqref{sde} backwad Kolmogorov equation
\begin{equation}
\label{eqp1}
(\partial_t-\Delta + b \cdot \nabla)u=0
\end{equation}
on the torus, for the entire class of form-bounded vector fields $b \in \mathbf{F}_\delta$. 
Our main result (Theorem \ref{thm1}) is a well-posedness theory of \eqref{eqp1} in an Orlicz space 
situated between all $L^p$ and $L^\infty$. This space is essentially dictated by the drift term. 
Orlicz spaces are known to appear in analysis and the theory of PDEs in various borderline situations, 
cf.\,Trudinger's theorem, or see \cite{KM,M} regarding Orlicz spaces arising in the study of 
dynamics of compressible fluids. Theorem \ref{thm1} gives another examples of a borderline situation where an Orlicz space appears. We consider Theorem \ref{thm1} to be the first step towards a theory of \eqref{eqp1} for $\delta=4$ that should establish a link with singular SDEs.

\begin{remark}[On interacting particle systems] 
\label{mv_rem}
SDE \eqref{sde} considered in $\mathbb R^{Nd}$, with $B_t=(B_t^1,\dots,B_t^N)$ being the vector of $N$ independent $d$-dimensional Brownian motions, $N \geq 2$, $X_t=(Y_t^1,\dots,Y_t^N)$ where $Y^i_t$ is the position of the $i$-th particle in $\mathbb R^d$ at time $t$, and drift $b=(b^1,\dots,b^N):\mathbb R^{Nd} \rightarrow \mathbb R^{Nd}$, defined by
\begin{equation}
\label{b_def}
b^i(y_1,\dots,y_N):=\frac{1}{N}\sum_{j=1,j  \neq i}^N K(y_i-y_j), \quad 1 \leq i \leq N, \quad y_j \in \mathbb R^d,
\end{equation}
describes the dynamics of $N$ interacting particles with the interaction kernel $K:\mathbb R^d \rightarrow \mathbb R^d$. If $K \in \mathbf{F}_\kappa$ (on $\mathbb R^d$), then a simple calculation shows that $b$ defined by \eqref{b_def} satisfies
$$
b \in \mathbf{F}_\delta \text{ on }\mathbb R^{Nd}, \quad \delta=4\biggl(\frac{N}{N-1}\biggr)^2,
$$
moreover, for $K(y)=\pm \sqrt{\delta}\frac{d-2}{2}|y|^{-2}y$, $y \in \mathbb R^d$, one can improve (i.e.\,lower) the form-bound $\delta$ using the multi-particle Hardy inequality of \cite{HHLT}. See \cite{Ki_particle}. Theorem \ref{thm1} below thus applies to the backward Kolmogorov equation that is behind the multi-particle stochastic system \eqref{sde}, \eqref{b_def}. 

As $N \rightarrow \infty$, the empirical measures of the trajectories of $Y^i_t$ are expected to converge (cf.\,\cite{FJ}) to solution of the McKean-Vlasov SDE whose law $\rho_t(\cdot)$ satisfies the non-linear McKean-Vlasov PDE in $\mathbb R^d$
\begin{equation}
\label{mv}
\partial_t \rho - \Delta \rho - \nabla \cdot (\tilde{b}\rho)=0, \quad \tilde{b}(t):=K \ast \rho(t).
\end{equation}
with the initial probability distribution $\rho_0 \geq 0$, $\int_{\mathbb R^d }\rho_0(y)dy=1$. This $b$ is itself a time-dependent form-bounded vector field, i.e.\,for every $t \geq 0$ (let $\langle\,\rangle$ denote the integration over $\mathbb R^d$)
\begin{align*}
\langle |\tilde{b}(t)|^2\varphi^2\rangle & = \big\langle |\langle K(\cdot-z)\rho_t(z)\rangle_z|^2\varphi^2\big\rangle  \\
& (\text{apply Cauchy-Schwartz inequality in $z$ and use $\langle \rho_t(z)\rangle_z=1$}) \\
& \leq \big\langle \langle |K(\cdot-z)|^2\rho_t(z)\rangle_z\varphi^2\big\rangle  = \big\langle \langle |K(\cdot-z)|^2 \varphi^2\rangle \rho_t(z)\big\rangle_z \\
& (\text{apply $K \in \mathbf{F}_\kappa(\mathbb R^d)$ and use again $\langle \rho_t(z)\rangle_z=1$}) \\
& \leq \kappa\langle |\nabla \varphi|^2\rangle + c_\kappa \langle \varphi^2 \rangle.
\end{align*}
(At this point we should add that our Theorem \ref{thm1} extends easily to time-dependent form-bounded $b$, with propagators instead of semigroups, see comments in Section \ref{rem_sect}.) We exploit this in Section \ref{blow_sect} to rule out using Theorem \ref{thm1} the blow up of solution of the forward Kolmogorov equation of type \eqref{mv}, i.e.\,the formation of a delta-function in $\rho_t$ in finite time, under appropriate condition on the initial probability distribution. The questions of well-posedness and blow up of solutions of \eqref{mv} received a lot of attention in the past few years, see \cite{CP,FJ} and recent papers \cite{FT,T} dealing with the borderline strengths of attraction. Our results, however, are parallel to these developments, at least at the moment: there the authors work in dimension $d=2$ with the Hardy-type attracting interaction $K(y)=c|y|^{-2}y$, $y \in \mathbb R^d$, as needed for applications to the famous Keller-Segel model of chemotaxis. Moreover, their proofs rely on the special structure of $c|y|^{-2}y$.  We, on the other hand, deal with the entire class of form-bounded vector fields in dimensions $d \geq 3$. (Speaking of the critical strengths of interactions, see also \cite{BJW} who proved quantitative propagation of chaos for \eqref{mv} for a class of singular vector fields of gradient form in dimensions $d \geq 2$, covering the entire range, up to the equality, of the strengths of attracting interactions.)
\end{remark}

\begin{remark}[On the sub-critical condition $\delta<4$]
\label{subcr_rem}

Let us explain where does condition $\delta<4$ come from.
The authors of \cite{KS} proved, among many other results, that one can construct a strongly continuous semigroup corresponding to parabolic equation \eqref{eqp1} for all $\delta<4$ by working in $L^p$, $p > \frac{2}{2-\sqrt{\delta}}$. The following calculation illustrates this. Consider initial-value problem
\begin{equation}
\label{cauchy0}
\left\{
\begin{array}{rr}
(\partial_t - \Delta +b \cdot \nabla)u=0 \text{ on } [0,\infty[ \times \mathbb R^d, \\[2mm]
 u(0,\cdot)=f(\cdot),
\end{array}
\right.
\end{equation} 
where $b$ and $f$ are assumed to be smooth (but the constants in the estimates should not depend on the smoothness of  $b$ and $f$). Replacing $u$ by $v:=u e^{-\lambda t}$, $\lambda \geq 0$, we can deal with initial-value problem $(\lambda + \partial_t - \Delta +b \cdot \nabla)v=0$, $v(0)=f.$ Multiply this equation by $v^{p-1}$, where, without loss of generality, $p$ is rational with odd denominator (so that we can raise negative functions to power $p-1$), and integrate by parts: 
$$
\lambda \langle v^p\rangle + \frac{1}{p}\langle \partial_t v^p\rangle + \frac{4(p-1)}{p^2}\langle |\nabla v^{\frac{p}{2}}|^2 \rangle + \frac{2}{p}\langle b \cdot \nabla v^{\frac{p}{2}},v^{\frac{p}{2}}\rangle=0.
$$
Applying quadratic inequality in the last term (and multiplying by $p$), we arrive at
$$
p\lambda \langle v^p\rangle + \langle \partial_t v^p\rangle + \frac{4(p-1)}{p}\langle |\nabla v^{\frac{p}{2}}|^2 \rangle \leq \alpha \langle |b|^2,v^p \rangle + \frac{1}{\alpha} \langle |\nabla v^{\frac{p}{2}}|^2 \rangle
$$
Now, applying $b \in \mathbf{F}_\delta$ and selecting $\alpha=\frac{1}{\sqrt{\delta}}$, we obtain the energy inequality
\begin{equation}
\label{ei}
\biggl[p\lambda  - \frac{c_\delta}{\sqrt{\delta}}\biggr]\langle v^p\rangle + \langle \partial_t v^p\rangle + \biggl[\frac{4(p-1)}{p}-2\sqrt{\delta} \biggr]\langle |\nabla v^{\frac{p}{2}}|^2 \rangle \leq 0, \quad \lambda \geq \frac{c_\delta}{p\sqrt{\delta}}.
\end{equation}
Thus, in order keep the dispersion term one needs $\frac{4(p-1)}{p}-2\sqrt{\delta}> 0$, i.e.\,$p > \frac{2}{2-\sqrt{\delta}}$, hence the requirement $\delta<4$. (And even if one is ready to sacrifice the dispersion term, in order to control $\langle v^p\rangle$ one still needs it to be zero, hence one needs $\delta<4$.) 

Under the sub-critical condition $\delta<4$, one can remove the assumption of smoothness of $b$ and construct the corresponding strongly continuous quasi-contraction semigroup $e^{-t\Lambda_p(b)}f:=u(t)$ in $L^p$, see \cite{KS,KiS_super}.

Although the interval of quasi contraction solvability  $p \in ]\frac{2}{2-\sqrt{\delta}},\infty[$ tends to $\varnothing$ as $\delta \uparrow 4$, this interval can be extended to the interval of quasi-bounded solvability $q \in ]\frac{2}{2-\frac{d-2}{d}\sqrt{\delta}},\infty[$:
$$
\|e^{-t\Lambda_q(b)}\|_q \leq M_{q,\delta} e^{\lambda_{q,\delta} t}\|f\|_q, \quad t>0,
$$
see \cite{KiS_super} where the authors also demonstrated that the interval of quasi-bounded solvability is sharp. Notice that as $\delta \uparrow 4$ the interval of quasi-bounded solvability tends to $]\frac{d}{2},\infty[$. However,  $M_{q,\delta} \rightarrow \infty$ as $\delta \uparrow 4$, so this still does not give a strongly continuous semigroup for \eqref{eqp1} for $\delta=4$ in any $L^q$ with finite $q$. 

It should be added that one has a priori bound
$\|u(t)\|_\infty \leq \|f\|_\infty$ for solution $u$ of \eqref{cauchy0} regardless of the value of $\delta$.
However, the existing methods of constructing Feller semigroup for $-\Delta + b \cdot \nabla$ on $\mathbb R^d$ (i.e.\,a strongly continuous semigroup in the space of continuous functions vanishing at infinity endowed with norm $\|\cdot\|_\infty$) are based on the regularity theory of \eqref{cauchy0}, such as strong gradient bounds or De Giorgi's method, see \cite{KS,Ki,Ki_survey}. So, these methods require the dispersion term in \eqref{ei}, and hence the strict inequality $\delta<4$. (That said, in Theorem \ref{thm1} we still prove an energy inequality for $\delta=4$ with some weakened dispersion term.)

\end{remark}

\subsection{Main result. Critical case $\delta=4$}
\label{res_sect}
In the rest of the paper we work over  $d$-dimensional torus $\Pi^d$ obtained as the quotient of $[-\frac{1}{2},\frac{1}{2}]^d$. This is not a technical assumption since the volume of the torus enters the estimates; the case of $\mathbb R^d$ requires separate study. Still, since $\delta$ measures the magnitude of \textit{local} singularities of $b$, working on a torus is sufficient for the purposes of this paper.

The functions/vector fields on $\Pi^d$ are identified with $1$-periodic functions/vector fields on $\mathbb R^d$. Let $dx$ denote the Lebesgue measure on $\Pi^d$. Given a Borel measurable function $f:\Pi^d \rightarrow \mathbb R$, we put
$$
\langle f\rangle:=\int_{\Pi^d} f(x)dx,\qquad \langle f,g\rangle:=\langle fg \rangle.
$$
We have $|\Pi^d|=\langle 1 \rangle=1$. 
Let $\|\cdot\|_p$ denote the norm in $L^p\equiv L^p(\Pi^d,dx)$. Put $C^\infty:=C^\infty(\Pi^d)$.

We now introduce the ``critical'' Orlicz space on the torus.
Define the gauge function
\begin{align*}
\Phi(y)&:=\cosh(y)-1, \quad y \in \mathbb R
\end{align*}
It is clear that $\Phi(y)=\Phi(|y|)$ is convex on $\mathbb R$, $\Phi(y)=0$ if and only if $y=0$, $\Phi(y)/y \rightarrow 0$ if $y \rightarrow 0$ and $\Phi(y)/y \rightarrow \infty$ if $y \rightarrow \infty$. Therefore,
\begin{equation}
\label{Phi_norm}
\|f\|_\Phi:=\inf\big\{c>0 \mid \langle \Phi\bigl(\frac{f}{c}\bigr) \rangle \leq 1\big\}<\infty,
\end{equation}
is a norm and we can define $L_\Phi$ to be the closure of $C^\infty(\Pi^d)$ in the Orlicz space $\mathcal L_\Phi$, see \cite[Ch.\,8]{AF}.
With some abuse of terminology, we will call $L_\Phi$ itself the Orlicz space with gauge function $\Phi$.

 It follows from the Taylor series representation 
$\Phi(y) =\sum_{k=1}^\infty \frac{y^{2k}}{(2k)!}$
 that 
\begin{equation}
\label{Lpnorm}
\|\cdot\|_\Phi \geq \frac{1}{(2p)!}\|\cdot\|_{2p}, \quad p=1,2,\dots,
\end{equation}
i.e.\,$\|\cdot\|_\Phi$ is stronger than $\|\cdot\|_p$ regardless of how large (finite) $p$ is. It is, however, weaker than $\|\cdot\|_\infty$.

The definition of form-bounded vector fields on the torus does not change.

\begin{definition*}
A vector field $b \in [L^2(\Pi^d)]^d$ is called form-bounded if  there exists constant $\delta>0$ such that quadratic form inequality 
$$
\|b\varphi\|_2^2 \leq \delta\|\nabla \varphi\|_2^2+c_\delta\|\varphi\|_2^2 \quad \forall \varphi \in C^\infty
$$
holds
for some constant $c_\delta$. This will also be written, with some abuse of notation, as $b \in \mathbf{F}_\delta$.
\end{definition*}

The examples of form-bounded vector fields on $\mathbb R^d$ remain essentially unchanged when one transitions to $\Pi^d$. For relevant papers, we refer to \cite{BO,G}. One necessary comment, however, is related to the fact that  when working over the torus the Hardy inequality is often obtained for functions having zero average. This does not cause a problem for us (also because, in principle, we could also work with zero average functions in Theorem \ref{thm1}). For instance, a Hardy inequality in \cite{G}
$$
\int_{\Pi^d}V(x)|\varphi(x)-\langle\varphi\rangle|^2dx \leq \frac{3d^2+8d+4}{d(d-2)^2}\int_{\Pi^d}|\nabla \varphi(x)|^2 dx, \quad V(x):=\frac{1}{\sum_{j=1}^d \sin^2(x_j/2)},
$$
where $\langle \varphi \rangle:=\int_{\Pi^d}\varphi(x) dx$ coincides with the average of $\varphi$ on $\Pi^d$,
gives us, using Cauchy-Schwartz inequality,
$$
\int_{\Pi^d} V(x)\varphi^2(x) \leq \frac{3d^2+8d+4}{d(d-2)^2}\int_{\Pi^d}|\nabla \varphi(x)|^2 dx + c\int_{\Pi^d}\varphi^2(x) dx
$$
for all $\varphi \in C^\infty(\Pi^d)$,
where constant $c=C(1+\int_{\Pi^d}V^2(y)dy)$ with $\int_{\Pi^d}V^2(y)dy<\infty$ in dimensions $d \geq 5$.

\medskip

Given a $b \in \mathbf{F}_\delta$, we put $$b_n:=E_{\varepsilon_n}b,$$ where $\{\varepsilon_n\} \downarrow 0$ is fixed arbitrarily and $E_\varepsilon:=e^{\varepsilon \Delta}$  is the De Giorgi mollifier on $\Pi^d$.
Then 
\begin{equation}
\label{conv0}
b_n \rightarrow b \quad \text{ in } L^2(\Pi^d)
\end{equation}
 and, importantly, $b_n \in \mathbf{F}_\delta$ with the same constant $c_\delta$ as the one for $b$, see the beginning of the proof of Theorem \ref{thm1}.

By the standard theory, for every $ n=1,2,\dots$, for every $f \in C^\infty(\Pi^d)$, 
Cauchy problem for the backward Kolmogorov equation
$$
(\partial_t - \Delta + b_n\cdot \nabla)u_n=0, \quad u_n(0)=f
$$
has unique classical solution, moreover, the operators $e^{-t\Lambda(b_n)}$ defined by $$e^{-t\Lambda(b_n)}f:=u_n(t), \quad t \geq 0,$$ constitute a semigroup that is strongly continuous on smooth $f$ in the norm of $L_\Phi$ since it is strongly continuous in the stronger norm of $L^\infty$.

\begin{theorem} 
\label{thm1}
Let $b \in \mathbf{F}_\delta=\mathbf{F}_\delta(\Pi^d)$, $0<\delta \leq 4$. The following are true:

\smallskip

{\rm (\textit{i})} There exists a strongly continuous quasi contraction Markov semigroup $e^{-t\Lambda(b)}$ on $L_\Phi$, 
$$
\|e^{-t\Lambda(b)}f\|_\Phi \leq e^{2\frac{c_\delta}{\sqrt{\delta}} t}\|f\|_\Phi, \quad t \geq 0,
$$
such that, for every $f \in C^\infty$, 
\begin{equation}
\label{conv5}
\|e^{-t\Lambda(b)}f-e^{-t\Lambda(b_n)}f\|_{\Phi} \rightarrow 0 \quad \text{ as } n \rightarrow \infty \text{ loc.\,uniformly in $t \geq 0$}.
\end{equation}
The generator $\Lambda(b)$ of the semigroup is the appropriate operator realization of the formal operator $-\Delta + b \cdot \nabla$ in $L_\Phi$, so $u(t):=e^{-t\Lambda(b)}f$, $f \in L_\Phi$ is the strong solution of the Cauchy problem for the backward Kolmogorov equation $(\partial_t-\Delta + b \cdot \nabla)u=0$, $u(0)=f$ in $L_\Phi$.

\smallskip

{\rm (\textit{ii})} This semigroup is  unique in the sense that for any other the choice of smooth vector fields $b_n \rightarrow b$ in $L^2$ that do not increase constants $\delta$, $c_\delta$ we have convergence \eqref{conv5} to the same limit semigroup $e^{-t\Lambda(b_n)}$.

\smallskip

\smallskip

{\rm (\textit{iii})} The function $v:=(\mu+\Lambda(b))^{-1}f$, $f \in C^\infty(\Pi^d)$, $\mu > \frac{2c_\delta}{\sqrt{\delta}}$ is the unique weak solution to elliptic Kolmogorov equation
\begin{equation}
\label{eq_el3}
(\mu-\Delta + b \cdot \nabla)v=f,
\end{equation}
i.e.\,$v \in W^{1,2}(\Pi^{d}) \cap L^\infty(\Pi^{d})$ and
$$
\mu \langle v,\varphi \rangle + \langle \nabla v,\nabla \varphi\rangle +  \langle b \cdot \nabla v,\varphi\rangle=\langle f,\varphi\rangle
$$
for every $\varphi \in W^{1,2}(\Pi^{d}) \cap L^\infty(\Pi^{d})$. 

\smallskip

{\rm (\textit{iv})} Let $p \geq 2$ be rational with odd denominator. The following a priori energy inequality holds for $u=e^{-t\Lambda(b_n)}f$:
\begin{align*}
\sup_{s \in [0,t]}\langle e^{u^{p}(s)}  \rangle +  4\frac{(p-1)}{p}\int_0^t \langle (\nabla u^{\frac{p}{2}})^2e^{u^{p}}\rangle ds & + 2(2-\sqrt{\delta}) \int_0^t\langle (\nabla e^{\frac{u^{p}}{2}})^2 \rangle ds \\
&\leq \langle e^{f^{p}}  \rangle + \frac{c_\delta}{\sqrt{\delta}}\int_0^t \langle e^{u^{p}} \rangle ds.
\end{align*}
In particular, for $\delta=4$,
$$
 \frac{1}{2}\sup_{s \in [0,t]}\langle e^{u^{p}(s)}  \rangle +  4\frac{(p-1)}{p}\int_0^t \langle (\nabla u^{\frac{p}{2}})^2e^{u^{p}}\rangle ds \leq \langle e^{f^{p}}  \rangle, \quad p=2,4,\dots
$$
provided $\frac{c_\delta}{\sqrt{\delta}}t<\frac{1}{2}$; the last constraint can be removed using the semigroup property.

\end{theorem}

It is not yet clear what kind of well-posedness result is valid for SDE \eqref{sde} in the critical case $\delta=4$.
Of course, in the sub-critical case $\delta<4$ a lot more is known. Other than the result from \cite{KiS} on the weak solvability of SDE \eqref{sde} discussed in the beginning of the paper, let us mention \cite{KiS_MAAN} dealing with similar thresholds in the proofs of Gaussian lower and/or upper heat kernel bounds for \eqref{eqp1} and the general divergence-form equation \eqref{eqp2}.

In Theorem \ref{thm1} we could work over any compact Riemannian manifold without boundary. However, as was mentioned above, the finiteness of the volume of the manifold is important since the volume enters the estimates.

The last assertion of Theorem \ref{thm1} is noteworthy: at the first sight, it seems like the possibility to pass to $\delta=4$ comes at the cost of killing off the dispersion term. Nevertheless, it turns out that some gradient estimates persist even for $\delta=4$.

A crucial feature of Theorem \ref{thm1} is that it covers the entire class of form-bounded vector fields for the critical value of $\delta$ and not just some of its representatives as e.g.\,Hardy drift \eqref{hardy}.
In fact, \eqref{hardy} is  better than a typical representative of $\mathbf{F}_\delta$ since, on a bounded domain, such $b$ satisfies an ``improved form-boundedness condition''
$$
c\|\varphi\|^2_{2j}+\|b\varphi\|_2^2 \leq \delta\|\nabla \varphi\|_2^2 + c_\delta\|\varphi\|_2^2, \quad j<\frac{d}{d-2}, c>0.
$$
The latter is a re-statement of the improved Hardy inequality due to \cite{BV} (there $c_\delta=0$). 
Also, for this $b$, the corresponding forward Kolmogorov operator admits, at least formally, an explicit invariant measure, which opens up other ways for studying this equation; see \cite{BKRS} in this regard.

The class of form-bounded vector fields appears naturally in connection with divergence-form parabolic  equation
\begin{equation}
\label{eqp2}
(\partial_t - \nabla \cdot a \cdot \nabla + b \cdot \nabla)u=0,
\end{equation}
where $a$ is a symmetric uniformly elliptic matrix, i.e.\,$\sigma I \leq a \in [L^\infty]^{d \times d}$ for some $\sigma>0$.  More specifically, $b \in \mathbf{F}_\delta$ with $\delta<\sigma $  allows to verify coercivity of the corresponding sesquilinear form on $W^{1,2}$ and hence to apply the Kato-Lions-Lax-Milgram-Nelson (KLMN) theorem \cite[Ch.VI]{K}, which yields the strong solution theory  of $-\nabla \cdot a \cdot \nabla + b \cdot \nabla$ in $L^2$. Moreover, if one focuses on the assumptions on $b$ in terms of $|b|$ only, as we do in the present paper, then condition $b \in \mathbf{F}_\delta$ is in fact necessary for such $L^2$ theory to exist, see \cite{MV}.

\subsection{Forward Kolmogorov equation and absence of blow up} 
\label{blow_sect}
We discuss briefly the forward Kolmogorov equation
\begin{equation}
\label{f_eq}
\partial_t \rho_t - \Delta \rho_t - \nabla \cdot (b\rho_t)=0, \quad b \in \mathbf{F}_\delta, \quad \delta \leq 4,
\end{equation}
with the initial condition $\rho_0=g$, $0 \leq g$, $\langle g\rangle=1$.
The method of constructing semigroup $e^{-t\Lambda(b)}$ in the  proof of Theorem \ref{thm1} (i.e.\,verifying Cauchy criterion in $L^\infty([0,1],L_\Phi)$) does two things at the same time: it produces, for each $t>0$, the operators $e^{-t\Lambda} \in \mathcal B(L_\Phi)$, and yields the strong continuity of the semigroup in $L_\Phi$. It is not clear how to adapt this argument to \eqref{f_eq}. In fact, the strong continuity of the semigroup seems to be an obstacle here.

In view of Theorem \ref{thm1}, one expects \eqref{f_eq} to be posed in the dual of the Orlicz space $L_\Phi$. 
By a classical result (see \cite[Ch.\,8]{AF}), the dual of $L_\Phi$ is isomorphic and homeomorphic to $L_{\Psi}$, where $\Psi$ is the complimentary function (of logarithmic type) of $\Phi=\cosh - 1$. By the Phillips theorem, the adjoint operators $(e^{-t\Lambda})^\ast \in \mathcal B(L_\Psi)$ constitute a semigroup in $L_\Psi$, although this semigroup does not have to be strongly continuous, see \cite[Sect.\,IX.13]{Yos}. We thus can pose, even if formally, 
$$
\rho_t:=(e^{-t\Lambda})^\ast g, \quad \text{ for all }t>0
$$
provided that the initial probability distribution $g$ satisfies additional condition $\|g\|_{L_\Psi}<\infty$.
Returning to our discussion in Remark \ref{mv_rem}, we note that since the delta-function does not belong to $L_\Psi$, condition $\|g\|_{L_\Psi}<\infty$ excludes the blow up of $\rho_t$ in finite time.

\subsection{Further remarks}

\label{rem_sect}

1.~Speaking of the paper \cite{KS} that introduced strong solution theory of \eqref{eqp1} with $\delta<4$, the goal of the authors there was to detect the dependence of the regularity properties of solutions of \eqref{eqp1} on the value of $\delta$, which they did by showing that the strongly continuous semigroup for \eqref{eqp1} exists in $L^p$ for $p>\frac{2}{2-\sqrt{\delta}}$. But to reach $\delta=4$ one needs to work in a space that ``does not sense'' $0<\delta \leq 4$, such as the Orlicz space $L_\Phi$, $\Phi=\cosh - 1$.

\smallskip

2.~The proof of Theorem \ref{thm1} also works for form-bounded $b=b(t,x)$, i.e.\,$b \in L^2([0,\infty[ \times \Pi^d)$ and for a.e.\,$t \geq 0$
$$
\|b(t)\varphi\|_2^2 \leq \delta\|\nabla \varphi\|_2^2+g_\delta(t)\|\varphi\|_2^2 \quad \forall \varphi \in C^\infty
$$
for a function $0 \leq g_\delta \in L^1_{\loc}[0,\infty[$.

\smallskip

3.~Estimate (\textit{iv}) of Theorem \ref{thm1} extends right away to solutions of divergence form equation
 \eqref{eqp2} where the uniformly elliptic matrix $a \geq \sigma I$ and the vector field $b$ are assumed to be smooth and bounded, but the constants in  (\textit{iv}) do not depend on the smoothness of $a$, $b$, or the boundedness of $b$ (i.e.\,these are a priori estimates).

Define vector field $\nabla a$ by $(\nabla a)_k:=\sum_{i=1}^d \nabla_i a_{ik}$ ($1 \leq k \leq d$).
One obtains similar results for solutions of the non-divergence form equation via the representation
$$
- a \cdot \nabla^2 + b \cdot \nabla := - \nabla \cdot a \cdot \nabla + (\nabla a+b) \cdot \nabla,
$$
provided that $\nabla a+ b \in \mathbf{F}_\delta$, $ \delta \leq 4\sigma^2.$ This condition includes some matrices $a$ having critical discontinuities. For instance, working for simplicity on $\mathbb R^d$, if $a=I+c x \otimes x|x|^{-2}$, $c>-1$, then $\nabla a=c(d-1)|x|^{-2}x \in \mathbf{F}_{\delta}$ with $\delta=4c^2(d-1)^2/(d-2)^2$, see \cite[Sect.\,4]{KiS_super}. Since we can now reach the critical value of the form-bound $\delta$, we can treat previously inaccessible, at least in the presence of a form-bounded drift, magnitude of critical discontinuities of $a$.

\bigskip

\section{Proof of Theorem \ref{thm1}}

Before we begin the proof, let us first show that $b_n:=E_{\varepsilon_n}b$ (De Giorgi mollifier on $\Pi^d$) indeed do not increase constants $\delta$ and $c_\delta$ of $b$. Arguing as in \cite{KiS_MAAN}, we note that
$|b_\varepsilon|\leq\sqrt{E_\varepsilon |b|^2}$ and so
\begin{align*}
\|b_\varepsilon \varphi\|_2^2 &\leq \langle E_\varepsilon |b|^2,\varphi^2\rangle=\|b\sqrt{E_\varepsilon \varphi^2}\|^2_2 \leq \delta\|\nabla\sqrt{E_\varepsilon \varphi ^2}\|_2^2+c_\delta\|\varphi\|_2^2,
\end{align*}
where
\begin{align*}
\|\nabla\sqrt{E_\varepsilon|\varphi|^2}\|_2 & =\big\|\frac{E_\varepsilon(|\varphi||\nabla|\varphi|)}{\sqrt{E_\varepsilon|\varphi|^2}}\big\|_2\\
&\leq \|\sqrt{E_\varepsilon|\nabla |\varphi||^2}\|_2=\|E_\varepsilon|\nabla |\varphi||^2\|_1^\frac{1}{2}\\
& \leq\|\nabla|f|\|_2\leq \|\nabla \varphi\|_2,
\end{align*}
i.e. $b_\varepsilon\in\mathbf{F}_{\delta}$ with the same $c_\delta$ as $b$, as claimed.

\medskip

Now, we start the proof of Theorem \ref{thm1}.
We replace $u_n$ by $v_n=e^{-\lambda t}u_n$, $\lambda=\frac{c_\delta}{\sqrt{\delta}}$, which satisfies
\begin{equation}
\label{eq1}
\left\{
\begin{array}{rr}
(\lambda + \partial_t - \Delta +b_n \cdot \nabla)v_n=0 \text{ on } [0,\infty[ \times \Pi^d, \\[2mm]
 v_n(0,\cdot)=f(\cdot) \in C^\infty(\Pi^d).
\end{array}
\right.
\end{equation}

(\textit{i}) 
Fix $n$ and write for brevity $v=v_n$. Let us first prove bound
\begin{equation}
\label{v_bd}
\|v(t)\|_\Phi \leq e^{\frac{c_\delta}{\sqrt{\delta}}t}\|f\|_\Phi.
\end{equation}
To this end, 
we multiply the equation by $e^v$ and integrate:
$$
\lambda \langle v,e^v\rangle + \langle \partial_t (e^v-1)\rangle + 4\langle (\nabla e^{\frac{v}{2}})^2\rangle + 2\langle b e^{\frac{v}{2}},\nabla e^{\frac{v}{2}}\rangle=0.
$$
By quadratic inequality,
\begin{equation}
\label{e1}
\lambda \langle v,e^v\rangle +  \langle \partial_t (e^v-1)\rangle + 4\langle (\nabla e^{\frac{v}{2}})^2\rangle \leq  \alpha \langle b^2 e^{v}\rangle + \frac{1}{\alpha}\langle (\nabla e^{\frac{v}{2}})^2\rangle.
\end{equation}
Applying $b \in \mathbf{F}_\delta$ and selecting $\alpha=\frac{1}{\sqrt{\delta}}$, we arrive at
\begin{equation}
\label{e2}
\lambda \langle v,e^v\rangle + \langle \partial_t (e^v-1)\rangle + (4-2\sqrt{\delta})\langle (\nabla e^{\frac{v}{2}})^2\rangle \leq \frac{c_\delta}{\sqrt{\delta}} \langle e^v\rangle.
\end{equation}
Using $\delta \leq 4$ (we are interested above all in $\delta=4$), one obtains, after integrating in time from $0$ to $t$:
\begin{equation*}
\lambda \int_0^t \langle v,e^{v}\rangle ds + \langle e^{v(t)}-1 \rangle \leq \langle e^{f}-1 \rangle + \frac{c_\delta}{\sqrt{\delta}} \int_0^t \langle e^{v}\rangle ds. 
\end{equation*}
Replacing in the last inequality $v$ by $-v$ and adding up the resulting inequalities, we obtain
\begin{equation*}
\lambda  \int_0^t \langle v\sinh(v) \rangle ds + \langle \cosh (v(t))-1\rangle \leq \langle \cosh (f)-1\rangle + \frac{c_\delta}{\sqrt{\delta}} \int_0^t \langle \cosh(v)\rangle ds.
\end{equation*}
Applying  $v\sinh(v) \geq \cosh(v)-1$, we arrive at 
\begin{equation}
\label{e3}
(\lambda-\frac{c_\delta}{\sqrt{\delta}}) \int_0^t \langle \cosh(v)-1 \rangle ds + \langle \cosh (v(t))-1\rangle \leq \langle \cosh (f)-1\rangle + \frac{c_\delta}{\sqrt{\delta}} t,
\end{equation}
where at the last step we have used the fact that volume $|\Pi^d|=1$. Recall that $\lambda= \frac{c_\delta}{\sqrt{\delta}}$.

Since our equation is linear, replacing everywhere $v$ by $\frac{v}{c}$, $c>0$, we have 
\begin{equation*}
(\lambda-\frac{c_\delta}{\sqrt{\delta}}) \int_0^t \langle \cosh(\frac{v}{c})-1 \rangle ds + \langle \cosh (\frac{v(t)}{c})-1\rangle \leq \langle \cosh (\frac{f}{c})-1\rangle + \frac{c_\delta}{\sqrt{\delta}} t.
\end{equation*}
Recalling our choice of $\lambda$, we have
\begin{equation*}
\langle \cosh (\frac{v(t)}{c})-1\rangle \leq \langle \cosh (\frac{f}{c})-1\rangle + \frac{c_\delta}{\sqrt{\delta}} t.
\end{equation*}
Let us fix $t$ and divide interval $[0,t]$ into $k$ subintervals: $[0,\frac{t}{k}]$, $[\frac{t}{k},\frac{2t}{k}],\dots$, \dots, $[\frac{(k-1)t}{n},t]$, where 
$k$ is large, i.e.\,is so that 
$$
\gamma:=\frac{c_\delta}{\sqrt{\delta}} \frac{t}{k}<1.
$$ 
Now, let $c_*>0$ be minimal such that $\langle \cosh (\frac{f}{(1-\gamma)c_*})-1\rangle=1$ (i.e.\,$\|f\|_\Phi=(1-\gamma)c_*$).
 Using the Taylor series expansion for $\cosh-1$, one sees that
$$
\cosh (\frac{f}{(1-\gamma)c_*})-1 \geq \frac{1}{1-\gamma}\biggl[\cosh (\frac{f}{c_*}))-1\biggr].
$$
So, $\langle \cosh (\frac{f}{c_*}))-1 \rangle \leq 1-\gamma$.
Therefore, 
$$
\langle \cosh (\frac{v(\frac{t}{k})}{c_*})-1\rangle \leq 1,
$$
and so $$\|v(\frac{t}{k})\|_\Phi \leq c_* \equiv \frac{1}{1-\gamma}\|f\|_\Phi \equiv \frac{1}{1-\frac{c_\delta}{\sqrt{\delta}}\frac{t}{k}}\|f\|_\Phi.$$
By the semigroup property, 
$$\|v(t)\|_\Phi \leq  (1-\frac{c_\delta}{\sqrt{\delta}}\frac{t}{k})^{-k}\|f\|_\Phi.$$
Taking $k \rightarrow \infty$, we obtain $\|v(t)\|_\Phi \leq  e^{\frac{c_\delta}{\sqrt{\delta}}t}\|f\|_\Phi$, i.e.\,we have proved \eqref{v_bd}.

\medskip

Next, we prove the convergence result in (\textit{i}). It suffices to carry out the proof for solutions $\{v_n\}$ of \eqref{eq1}. In three steps:

\smallskip

Step 1.~First, let us note that $\nabla v_n$ are bounded in $L^2([0,1] \times \Pi^d)$ uniformly in $n$. Indeed, multiplying $(\lambda + \partial_t - \Delta + b_n \cdot \nabla )v_n=0$ by $v_n$ and integrating over $[0,t] \times \Pi^d$, $0<t \leq 1$, we obtain
$$
\lambda\int_0^t \langle v_n^2\rangle ds + \frac{1}{2}\langle v_n^2(t)\rangle - \frac{1}{2}\langle f^2 \rangle + \int_0^t \langle (\nabla v_n)^2\rangle ds =-\int_0^t \langle b_n \cdot \nabla v_n,v_n\rangle ds,
$$
$$
\frac{1}{2}\langle v^2_n(t)\rangle - \frac{1}{2}\langle f^2 \rangle + \int_0^t \langle (\nabla v_n)^2\rangle ds \leq \alpha \int_0^T \langle (\nabla v_n)^2\rangle ds + \frac{1}{4\alpha}\int_0^T \langle b_n^2 v_n^2\rangle ds,
$$
where, by $\|v_n(s)\|_\infty \leq \|f\|_\infty$, $s \in [0,T]$,
$$
\int_0^t \langle b_n^2 v_n^2\rangle ds \leq \sup_n\int_0^t \|b_n\|_2^2 ds \|f\|^2_\infty =: C_0 \|f\|^2_\infty
$$
(in view of \eqref{conv0}, $C_0<\infty$). Hence, selecting above e.g.\,$\alpha=\frac{1}{2}$, we obtain
$$
\frac{1}{2}\langle v^2_n(t)\rangle + \frac{1}{2} \int_0^t \langle (\nabla v_n)^2\rangle ds \leq \frac{1}{2}\langle f^2 \rangle + \frac{1}{2}C_0  \|f\|^2_\infty.
$$
In particular,
\begin{equation}
\label{est4}
\sup_n \int_0^t \langle (\nabla v_n)^2\rangle ds \leq \|f\|_2^2 + C_0 \|f\|^2_\infty.
\end{equation}
(At this step we actually do not need positive $\lambda$, but we will need it at the next step.)

\medskip

Step 2.~ Let us show that $v_n-v_m \rightarrow 0$ in $L_\Phi$ as $n,m \rightarrow \infty$ uniformly in $t \in [0,T]$, where $0 <T \leq 1$ will be chosen later. (At the next step we will define the sought semigroup on $[0,T]$ as the limit of $v_n$.)

Put $$h:=\frac{v_n-v_m}{c}, \quad c>0.$$ We have
\begin{equation}
\label{eq}
\lambda h + \partial_t h - \Delta h + b_n \cdot \nabla h + (b_n-b_m) \cdot c^{-1}\nabla v_m=0, \quad h(0)=0.
\end{equation}
We multiply by $e^h$ and integrate by parts. The terms $\lambda h+ \partial_t h - \Delta h + b_n \cdot \nabla h$ are handled as in the beginning of the proof of (\textit{i}) (but with initial condition $h(0)=0$):
\begin{align}
(\lambda-\frac{c_\delta}{\sqrt{\delta}})\int_0^t \langle e^{h}-1\rangle ds & + \langle e^{h(t)}-1 \rangle + (4-2\sqrt{\delta})\int_0^t \langle (\nabla e^{\frac{h}{2}})^2\rangle ds \notag \\
& \leq - \int_0^t \langle (b_n-b_m) \cdot c^{-1}\nabla v_m,e^h \rangle ds + \frac{c_\delta}{\sqrt{\delta}}t. \label{est5}
\end{align}
Using $\|e^{h(s)}\|_\infty \leq e^{2c^{-1}\|f\|_\infty}$, we estimate:
\begin{align*}
\bigg|\int_0^t \langle (b_n-b_m) \cdot c^{-1}\nabla v_m,e^h \rangle ds\bigg|& \leq 
\bigg(\int_0^t \|b_n-b_m\|^2_2 ds\bigg)^{\frac{1}{2}} c^{-1} \bigg(\int_0^t \|\nabla v_m\|^2_2 ds\bigg)^{\frac{1}{2}}
 e^{2c^{-1}\|f\|_\infty} \\
& (\text{use Step 1}) \\
& \leq \bigg(\int_0^t \|b_n-b_m\|^2_2 ds\bigg)^{\frac{1}{2}}c^{-1}\biggl( \|f\|_2^2 + C_0 \|f\|^2_\infty \biggr)^{\frac{1}{2}}e^{2c^{-1}\|f\|_\infty}.
\end{align*}
By \eqref{conv0}, $\int_0^t \|b_n-b_m\|^2_2 ds \rightarrow 0$ as $n,m \rightarrow \infty$.
So, for every $c>0$,
\begin{equation}
\label{conv}
\bigg|\int_0^t \langle (b_n-b_m) \cdot c^{-1}\nabla v_n,e^h \rangle ds\bigg| \rightarrow 0 \quad \text{ as }n,m \rightarrow \infty \text{ uniformly in $0 \leq t \leq T$}.
\end{equation}

Now, since $\delta \leq 4$, we have by \eqref{est5} (recall: $\lambda=\frac{c_\delta}{\sqrt{\delta}}$) and \eqref{conv}, for every fixed $c>0$, for all $\varepsilon>0$
$$
\sup_{t \in [0,T]}\langle e^{\frac{v_n(t)-v_m(t)}{c}}-1 \rangle \leq \varepsilon + \frac{c_\delta}{\sqrt{\delta}}T
$$
for all $n,m$ sufficiently large.

Repeating the previous argument for $-h$ and adding up the resulting inequalities, we obtain: for every fixed $c>0$, for all $\varepsilon>0$, 
$$
\sup_{t \in [0,T]}\langle \Phi \bigl(\frac{v_n(t)-v_m(t)}{c}) \rangle \leq \varepsilon + \frac{c_\delta}{\sqrt{\delta}}T
$$
for all $n,m$ sufficiently large. Selecting $T$ such that $\frac{c_\delta}{\sqrt{\delta}}T<1$, we thus obtain for every $c>0$, provided that $\varepsilon$ is chosen sufficiently small: $\sup_{t \in [0,T]}\langle \Phi \bigl(\frac{v_n(t)-v_m(t)}{c}) \rangle \leq 1$ for all $n,m$ large enough. Hence $\|v_n(t)-v_m(t)\|_\Phi \rightarrow 0$ as $n,m \rightarrow \infty$ uniformly in $0 \leq t \leq T$.

\medskip

Step 3.~Define $$S^t f:=L_\Phi\mbox{-}\lim_n e^{\frac{c_\delta}{\sqrt{\delta}} t}v_n(t) \equiv L_\Phi\mbox{-}\lim_n e^{-t\Lambda(b_n)}f, \quad t \in [0,T].$$
This is a continuous $L_\Phi$ valued function of $t \in [0,T]$. By passing to the limit in $n$ in $\|v_n(t)\|_\Phi \leq e^{\frac{c_\delta}{\sqrt{\delta}} t}\|f\|_\Phi$, see (\textit{i}), we obtain $\|S^t f\|_\Phi \leq e^{2\frac{c_\delta}{\sqrt{\delta}} t}\|f\|_\Phi$. The linearity of $S^t$ is evident. The semigroup property ($t,s \in [0,T]$): 
\begin{align*}
\|e^{-t\Lambda(b_n)}e^{-s\Lambda(b_n)}f - S^t S^s f\|_\Phi & \leq \|(e^{-t\Lambda(b_n)}(e^{-s \Lambda(b_n)}f-S^sf)\|_\Phi + \|(e^{-t \Lambda(b_n)}-S^t )S^s f)\|_\Phi \\
& \leq \|e^{-s \Lambda(b_n)}f-S^sf\|_\Phi + \|e^{-t \Lambda(b_n)}-S^t )S^s f\|_\Phi \rightarrow 0, \quad n \rightarrow \infty.
\end{align*}
On the other hand, $e^{-t\Lambda(b_n)}e^{-s\Lambda(b_n)}f=e^{-(t+s)\Lambda(b_n)}f \rightarrow S^{t+s}f$, and so the semigroup property follows.

We extend $S^t$ from $C^\infty$ to $L_\Phi$ via the standard density argument using $\|S^t f\|_\Phi \leq e^{2\frac{c_\delta}{\sqrt{\delta}} t}\|f\|_\Phi$.
Finally, we extend $S^t$ to all $t>0$ by postulating the semigroup property.

\smallskip

(\textit{ii}) This is clear from the construction of the semigroup via Cauchy's criterion. That is, in the proof of (\textit{ii}), say, we have two smooth approximations $\{b_n\}$, $\{b_n'\}$ of $b$ in $L^2$ not increasing the constants $\delta$, $c_\delta$ of $b$, such that, for a fixed initial function $f$, the corresponding solutions $v_n$, $v_n'$ converge to different limits. However, mixing $\{b_n\}$, $\{b_n'\}$, we obtain that the corresponding sequence of solutions is again a Cauchy sequence, and so the limits of $v_n$, $v_n'$ must coincide.

\smallskip

We prove (\textit{iii}) below.

\smallskip

(\textit{iv}) We multiply equation $(\partial_t - \Delta + b_n \cdot \nabla )u=0$ by $u^{p-1}e^{u^p}$ and integrate:
\begin{equation}
\label{eq10}
\frac{1}{p} \langle \partial_t e^{u^p}\rangle+ \langle (-\Delta u), u^{p-1}e^{u^p} \rangle + \langle b \cdot \nabla u, u^{p-1}e^{u^p}\rangle=0,
\end{equation}
where
\begin{align*}
\langle (-\Delta u), u^{p-1}e^{u^p} \rangle & = (p-1)\langle \nabla u, u^{p-2}(\nabla u)e^{u^p}\rangle + p\langle \nabla u,u^{p-1}e^{u^p} u^{p-1}\nabla u\rangle \\
& = \frac{4(p-1)}{p^2}\langle (\nabla u^{\frac{p}{2}})^2e^{u^p}\rangle + \frac{4}{p} \langle (\nabla e^{\frac{u^p}{2}})^2\rangle
\end{align*}
and
\begin{align}
\langle b \cdot \nabla u, u^{p-1}e^{u^p}\rangle & = \frac{2}{p}\langle b \cdot \nabla e^{\frac{u^p}{2}},e^{\frac{u^p}{2}}\rangle \notag \\
& \leq \frac{1}{p}\biggl(\alpha \langle |b|^2 e^{u^p} \rangle + \frac{1}{\alpha} \langle (\nabla e^{\frac{u^p}{2}})^2 \rangle \biggr) \notag \\
& \leq \frac{1}{p} \biggl(\alpha \delta + \frac{1}{\alpha} \biggr)\langle (\nabla e^{\frac{u^p}{2}})^2 \rangle + \frac{1}{p}\alpha c_\delta \langle e^{u^p} \rangle \qquad \alpha:=\frac{1}{\sqrt{\delta}} \notag \\
& \leq \frac{2}{p}\sqrt{\delta}\langle (\nabla e^{\frac{u^p}{2}})^2 \rangle + \frac{1}{p}\frac{c_\delta}{\sqrt{\delta}}\langle e^{u^p} \rangle. \label{a9}
\end{align}
Applying this in \eqref{eq10}, we obtain 
assertion (\textit{iv}).

\smallskip

(\textit{iii}) First, we show that $v=(\mu+\Lambda(b))^{-1}f$, $\mu>\frac{2c_\delta}{\sqrt{\delta}}$ where the resolvent $(\mu+\Lambda(b))^{-1}$ was constructed in (\textit{i}), is a weak solution to \eqref{eq_el3}. Let $v_n$ be the classical solution to $(\mu - \Delta+ b_n \cdot \nabla)v_n=f$. On the one hand, by assertion (\textit{i}), $v_n \rightarrow (\mu+\Lambda(b))^{-1}f$ in $L_\Phi$ as $n \rightarrow \infty$, and hence in $L^2$. On the other hand, the same argument as the one used in the proof of \eqref{est4} but applied to the last equation yields
$
\sup_n\|\nabla v_n\|_2<\infty
$
(one can also use the energy inequalities in (\textit{iv}), (\textit{v}), but this will do since we assume here that $f$ is bounded).
Therefore, we can extract a sub-sequence of $\{v_{n'}\}$ such that $\nabla v_{n'}$ converges weakly in $L^2$ and, since $\nabla$ is a closed operator, the limit is $\nabla v \in L^2$. We can now pass to the limit $n \rightarrow \infty$ in 
$$
\mu \langle v_n,\varphi \rangle + \langle \nabla v_n,\nabla \varphi\rangle + 
\langle b_n \cdot \nabla v_n,\varphi\rangle=\langle f,\varphi\rangle, \quad \varphi \in W^{1,2}(\Pi^{d}) \cap L^\infty(\Pi^d)
$$
using $b_n \rightarrow b$ in $L^2(\Pi^{d})$.

The uniqueness of weak solution $v$ now follows using a standard argument. Let $\tilde{v}$ be another weak solution of \eqref{eq_el3}. Set $w:=v-\tilde{v}$.
Then
$$
\mu \langle w,\varphi \rangle + \langle \nabla w,\nabla \varphi\rangle + 
\langle b \cdot \nabla w,\varphi\rangle=0, \quad \forall\,\varphi.
$$
Since $we^{w^2} \in W^{1,2}(\Pi^{d}) \cap L^\infty(\Pi^d)$, we can take $\varphi:=we^{w^2}$, so
$$
\mu \langle w^2e^{w^2}\rangle + \langle \nabla w,(\nabla w)e^{w^2}\rangle + 2\langle |\nabla e^{\frac{w^2}{2}}|^2\rangle=  - \langle  b \cdot \nabla w,we^{w^2}\rangle.
$$
Now, estimating the RHS as in \eqref{a9} (take $p=2$ there), we obtain
$$
\mu \langle w^2e^{w^2}\rangle  \leq \frac{c_\delta}{\sqrt{\delta}}\langle e^{w^2}\rangle.
$$
Out goal is to prove that $w=0$. To this end, we use the fact that $w$ satisfies a linear equation, while the last inequality is non-linear. That is, fix arbitrary $\varepsilon>0$, let $A_\varepsilon:=\{|w|>\varepsilon\}$. First, we note that, after replacing $\mu$ in the LHS by smaller constant $\frac{c_\delta}{\sqrt{\delta}}$ and then dividing both sides by it,
$$
 \langle \mathbf{1}_{A_\varepsilon}(w^2-1)e^{w^2}\rangle  \leq \langle \mathbf{1}_{A_\varepsilon^c}e^{w^2}\rangle.
$$
Next, we blow up the equation, i.e.\,replace $w$ by $cw$, $c>0$, obtaining
$$
 \langle \mathbf{1}_{A_\varepsilon}(c^2w^2-1)e^{c^2w^2}\rangle  \leq \langle \mathbf{1}_{A_\varepsilon^c}e^{c^2w^2}\rangle.
$$
In the LHS, we replace $w^2$ by smaller (on $A_\varepsilon$) value $\varepsilon^2$ (but not everywhere), and in the RHS we replace $w^2$ by larger (on $A^c_\varepsilon$) value $\varepsilon^2$, obtaining
$$
 \langle \mathbf{1}_{A_\varepsilon}(c^2\varepsilon^2-1)e^{c^2w^2}\rangle  \leq \langle \mathbf{1}_{A_\varepsilon^c}e^{c^2\varepsilon^2}\rangle.
$$
Now, let $c>\varepsilon^{-1}$, i.e.\,$c^2\varepsilon^2-1>0$, so that we can also replace in the LHS $e^{c^2w^2}$ by smaller $e^{c^2\varepsilon^2}$, arriving at
$$
 (c^2\varepsilon^2-1)e^{c^2\varepsilon^2}\langle \mathbf{1}_{A_\varepsilon}\rangle  \leq e^{c^2\varepsilon^2}\langle \mathbf{1}_{A_\varepsilon^c}\rangle.
$$
So, using $\langle \mathbf{1}_{A_\varepsilon^c}\rangle \leq 1$, we arrive at the inequality $ (c^2\varepsilon^2-1)\langle\mathbf{1}_{A_\varepsilon}\rangle \leq 1$ for all $\varepsilon^{-1}<c<\infty$, which is possible only if $\langle\mathbf{1}_{A_\varepsilon}\rangle=0$. Since $\varepsilon>0$ was arbitrary, $w=0$,
and the uniqueness of the weak solution follows.

\begin{remark}
We obtained a relatively simple proof of the existence and uniqueness of the weak solution because we were working in the elliptic setting and, importantly, restricted our attention to the bounded right-hand side $f$ of the elliptic equation. Arguably, this is rather sufficient for the probabilistic setting, but from the PDE of view it is desirable to extend assertion (\textit{iii}) to all $f \in L_\Phi$.
\end{remark}

\end{document}